\documentclass{amsart}
\usepackage{amssymb}
\usepackage{graphicx}
\def\({\left(}
\def\){\right)}

\newtheorem{lema}{Lemma}[section]

\newtheorem*{teorema*}{Theorem}

\newtheorem{remark}[lema]{Remark}

\newtheorem{lemma}{Lemma}[section]
\newtheorem{corollary}[lema]{Corollary}
\newtheorem{theorem}[lema]{Theorem}

\newtheorem{proposition}[lema]{Proposition}

\newtheorem{definition}[lema]{Definition}

  {\par \hfill \fbox{}}
  {\par \hfill \fbox{}}
  {\par \hfill }
  {\par \hfill }

\def\BB{\mathbb B}

\def\RR{\mathbb R}

\def\beq{\begin{equation}}
\def\eeq{\end{equation}}
\def\beginpf{\noindent{\bf Proof.} \quad}
\def\endpf{\qquad\hfill\rule{2.2mm}{2.2mm}\break}

\def\epsilon{\varepsilon}

\begin{document}

\title[COMPOSITION OPERATORS]{Slice regular Besov spaces of hyperholomorphic functions and composition operators}
\author{Sanjay Kumar}
\address{Department of Mathematics, Central University of Jammu,
Jammu 180 011, INDIA.} \email{sanjaykmath@gmail.com}
\author{Khalid Manzoor}
\address{Department of Mathematics, Central University of Jammu,
Jammu 180 011, INDIA.}
\email{khalidcuj14@gmail.com}


\date{\today.}

\begin{abstract}
In this paper we investigate  some basic  results on the slice regular Besov spaces  of hyperholomorphic functions   on the unit ball $\mathbb{B}.$ We  also characterize the boundedness, compactness and   find the essential norm estimates of composition operators between these spaces. 
\end{abstract}

\maketitle

\section{Introduction}
In the last ten years the theory of slice regular functions   is developed systemically in the papers
\cite{ colo10, gent13, gent12, gent121, gent11, gent111, gent1111, gent09, gent08, gent081, gent07,  gent06}. Slice
hyperholomorphic functions when defined and takes values in quaternions  are called slice regular, 
see \cite{ alpa14,  alpa15, chui71}. In case they are defined on the Euclidean space $\RR^{N+1}$ and 
takes values in the Clifford algebra $\RR_{\omega}$   they are called slice monogenic functions,
see \cite{birk29, birk36}.  Several function spaces  of the slice hyperholomorphic functions are 
studied. The quaternionic Hardy spaces are studied in  \cite{alpa15, alpa12, alpa13, arco15, arco14, arco, sarf13}. The 
Bergman spaces of slice hyperholomorphic functions are invesigated in \cite{colo13,  colo131, colo12}. For Fock spaces
in the slice hyperholomorphic settings,  see  \cite{alpa14}. Further,  weighted Bergman spaces, Bloch, Besov and Dirichlet 
spaces of slice hyperholomorphic functions on the unit ball are considered  in \cite{marco}.  D.~Alpay etc.al studied Schur
analysis in the slice hyperholomorphic setting see e.g., \cite{ abu15, abu151, alpa15, alpa13}  and references therein.
The study of slice hyperholomorphic functions have wide range of applications. For complete discussion of slice regular 
functions and their applications, we refer the book  \cite{gent1211} and a recent survey \cite{colo15}. \\
   For each  $q\in \mathbb{H},$ we can write
$ \mathbb{H}=\{q=x_0+ix_1+jx_2+kx_3,$ for all $x_1, x_2, x_3\in \mathbb{R}\},$  where $\{1, i, j, k\}$ form the  basis of  quaternions with   imaginary units
$i,j,k$  such that $i^2=j^2=k^2=-1, ij=-ji=k, jk=-kj=i,  ki=-ik=j.$   The Euclidean norm on  $ \mathbb{H}$ is given by  $|q| =\sqrt{q\bar q}=\sqrt{x^2_0+x^2_1+x^2_2+x^2_3},$  where  $\bar q=Re(q)-Im(q)=x_0-(ix_1+jx_2+kx_3)$ represents the conjugate of $q$ with $Re(q)=x_0,\;\; Im(q)=ix_1+jx_2+ kx_3$ and  the multiplicative inverse  $q^{-1}$ of  non-zero quarternion  $q$ is  given by $\displaystyle\frac{\bar q}{|q|^2}.$ An element $q$ in $ \mathbb{H}$ can be also written as  linear combination of  two complex numbers $q=(x_0+ix_1)+(x_2+ix_3)j.$
 By symbol  $\mathbb{S}$ we denote the two-dimensional  unit sphere of purely imaginary quaternions i.e,
$\mathbb{S}=\{q=ix_1+ jx_2 +kx_3 ~~\mbox{such that }~x_1^2+x_2^2 +x_3^2=1\}.$ If $I\in \mathbb{S}$ then $I^2=-1.$ Let $\mathbb{C}(i)$  be  the 
space generated by $ \{1,i\}.$  For any $i\in \mathbb{S}, $ let $\Omega_i= \Omega \cap \mathbb{C}_i,$ for some subset
(domain) $\Omega$ of $\mathbb{H}.$ For a nonreal quaternion $q$ we can write $q=x+Im(q)=x+I_q|Im(q)|=x+I_qy,$
where $x=x_0, \;\;y=ix_1+jx_2+kx_3$ with
$ I_q=\displaystyle\frac{Im(q)}{|Im(q)|},$ so it lies in the complex plane $\mathbb{C}(i).$ We define the slice regular functions on the open ball   $\mathbb{B}(0,1)=\mathbb{B}=\{q\in \mathbb{H}: |q|<1\}$ in $\mathbb{H}$ and $\mathbb{B}\cap \mathbb{C}(i)=\mathbb{B}_i$ denote the unit disk in the complex plane, for $i\in \mathbb{S}.$  The study of slice holomorphic functions is now an active area of research and lot of work is being done in this direction. For  slice holomorphic functions we refer to \cite{fcs2009, fcs13, gent07, marco} and  references therein.
Here, we collect some basic definitions and  basic results already obtained in the quaternionic-valued slice regular functions.

\begin{definition}
Let $\Omega$ be an open set  in $\mathbb{H}$.  A real differential   function $f:\Omega \to \mathbb{H}$ is said to be  (left) slice regular or slice hyperholomorphic  on ${\Omega}_i,$ if for every  $i\in \mathbb{S}$,  $$ \frac{\partial }{\partial x}f_i(x+iy)+i\frac{\partial }{\partial y}f_i(x+iy)=0,$$  where $f_i$ denote the restriction of $f$ to 
  $\Omega \cap \mathbb{C}(i).$ The  class of slice regular function on $\Omega$ is denoted  by $SR(\Omega).$
\end{definition}


\begin{lemma}\label{eq:1}
 \cite[Lemma 4.1.7 ]{fcs13}(Splitting Lemma) If $f$ is a slice regular function on the domain $\Omega,$ then for any 
$i,j\in  \mathbb{S},$  with $i\bot j$ there exists two holomorphic functions $f_1,f_2:\Omega_i \to\mathbb{C}(i) $ 
such that 
\beq \label{eq:2}
f_i(z)=f_1(z)+f_2(z)j; ~~\mbox{for any }~~z=x+iy \in {\Omega}_i.
\eeq
\end{lemma}
One of the most important property of the slice regular functions is their Representation Formula.  It only  holds on the open sets which are stated below.
\begin{definition}
 Let $\Omega$ be an open set in  $\mathbb{H}$. We say  $\Omega$ is Axially symmetric if for any  $q=x+I_qy \in \Omega,$ all the elements $x+iy$ is contained in $\Omega$, for all $i\in  \mathbb{S}$ and $\Omega$ is said to be slice domain if $\Omega\cap \mathbb{R}$ is non empty and  $\Omega\cap \mathbb{C}(i)$ is a domain in $\mathbb{C}(i)$ for all $i\in \mathbb{S}.$
\end{definition}
\begin{theorem} \label{eq:115}
  \cite[Theorem  4.3.2 ]{fcs13} (Representation Formula) Let  $f$ be a  slice regular function in the   symmetric slice  domain $\Omega\subseteq  \mathbb{H} $ and let $j\in  \mathbb{S}$. Then  for all $z=x+iy \in \Omega$ with $i\in  \mathbb{S},$ the following equality holds $$f(x+iy)=\frac{1}{2}\left\{(1-ij)f(x+jy)+(1+ij)f(x-jy)\right\}.$$
\end{theorem}
\begin{remark}
Let $i,j $ be orthogonal imaginary units in S and $\Omega$ be  an axillay symmetric slice domain. Then the Splitting Lemma and the  Representation formula  generate a class of  operators on the slice regular functions as follows:

$$Q_i:SR(\Omega)\to hol({\Omega}_i)+hol({\Omega}_i)j$$ $$Q_i: f\mapsto  f_1+f_2j$$

$$P_i: hol({\Omega}_i)+hol({\Omega}_i)j \to SR(\Omega)$$
$$P_i[f](q)=P_i[f](x+I_qy)=\frac{1}{2}[(1-I_qi)f(x+iy)+(1+I_qi)f(x-iy)].$$ Also, $$P_i\circ Q_i=\emph{I}_{SR(\Omega)}~ \mbox{and}~Q_i\circ P_i=\emph{I}_{SR( hol({\Omega}_i)+hol({\Omega}_i))},$$ where $\emph{I}$ is an identy operator.

\end{remark}

Since pointwise product of functions does not preserve slice regularity, a new multiplication operation for regular functions is defined. In the special case of power series, the regular product (or $\star-$product) of $f(q) = \sum_{n= 0}^{\infty}q^{n}a_{n}$ and  $g(q) = \sum_{n= 0}^{\infty}q^{n}b_{n}$ is $$ f \star g(q) =  \sum_{n \ge 0} q^{n} \sum_{k = 0}^{n} a_{k}b_{n-k}. $$
The $\star-$product is related to the standard pointwise product by the following formula.
\begin{theorem} 
  \cite[Proposition   2.4 ]{arco14}
Let $f, g$ be regular functions on $\BB.$ Then
$f \star g(q) = 0 $ if $f(q)  = 0 $ and $ f(q) g(f(q)^{-1} qf(q))$ if $ f(q) \neq 0.$ 
The reciprocal $f^{-\star}$ of a regular function $f (q) =\sum_{n= 0}^{\infty}q^{n}a_{n}$
with respect to the $\star-$product is
$$f^{\star}(q) = \frac{1}{f \star f^{c}(q)} f^{c}(q), $$ where $f^{c}(q) =\sum_{n= 0}^{\infty}q^{n}\overline{a_{n}} $ is the regular conjugate of $f.$ The function $f^{-\star}$
is regular on $\BB\setminus(q \in \BB | f \star f^{c} (q) = 0) $ and $f  \star f^{-\star} = 1$ there.
\end{theorem}
\section{Besov spaces}
\label{sec:2}
Now we define Besov space of slice hyperholomorphic functions on the unit ball $\mathbb{B}.$
 Let $\mathbb{D}$ be a unit disk in the complex plane $ \mathbb{C}$ and $dA$ be the normalized area measure on  $\mathbb{D}.$
 For $1<p< \infty,$ a holomorphic  function $f: \mathbb{D} \to  \mathbb{C}$ is said to be in  Besov space  $\mathfrak{B}_{p, \mathbb{C}}\mathbb{(D)}$ if $$\int_{\mathbb{D}}\left|(1-|z|^2)f'(z)\right|^p d\lambda(z) <\infty,$$
where $d\lambda(z) =\displaystyle\frac{dA(z)}{(1-|z|^2)^2}$ and is M\"{o}bius invariant measure on $\mathbb{D}$. The space $\mathfrak{B}_{p,  \mathbb{C}}$ is a Banach space under the norm $$\|f\|_{\mathfrak{B}_{p,  \mathbb{C}}}=|f(0)|+\left(\int_{\mathbb{D}}\left|(1-|z|^2)f'(z)\right|^p d\lambda(z)\right)^{\frac{1}{p}}. $$
For details on the Besov space of the unit disk one can refer to  \cite{arco02, zhu91, zhu2}  and references therein.
\begin{definition}
 Let $p>1$ and let $i\in \mathbb{S}$.  The quaternionic right linear space of slice regular functions $f$ is said to be the quaternionic slice  regular Besov  space on the unit ball $\mathbb{B},$ if
$$\sup_{i\in \mathbb{S}}\int_{\mathbb{B}_i}\left|(1-|q|^2)\frac{\partial f}{\partial x_0}(q)\right|^p d\lambda_i(q) <\infty,$$that is,
$$\mathfrak{B}_{p}=\{f\in SR(\mathbb{B}):\sup_{i\in \mathbb{S}}\int_{\mathbb{B}_i}\left|(1-|q|^2)\frac{\partial f}{\partial x_0}(q)\right|^p d\lambda_i(q) <\infty\},$$
where $d\lambda_i(q) =\displaystyle\frac{dA_i(q)}{(1-|q|^2)^2}$ and is M\"{o}bius invariant measure on $\mathbb{B}$. The space $\mathfrak{B}_{p}$ is a Banach space under the norm $$\|f\|_{\mathfrak{B}_{p}}=|f(0)|+\left(\sup_{i\in \mathbb{S}}\int_{\mathbb{B}_i}\left|(1-|q|^2)\frac{\partial f}{\partial x_0}(q)\right|^p d\lambda_i(q)\right)^{\frac{1}{p}} . $$
For details on Besov spaces of quaternions  holomorphic functions  one can refer  to  \cite{marco}.
 By space $\mathfrak{B}_{p,i}, p>1, $ we means the quaternionic right linear space of slice regular functions on the unit ball $\mathbb{B} $ such that $$\int_{\mathbb{B}_i}\left|(1-|z|^2) Q_i[f]'(z)\right|^p d\lambda_i(z) <\infty,$$ 
and the norm of this space is given by  $$\|f\|_{\mathfrak{B}_{p,i}}=|f(0)|+\left(\int_{\mathbb{B}_i}\left|(1-|z|^2) Q_i[f]'(z)\right|^p d\lambda_i(z)\right)^{\frac{1}{p}}$$
where $Q_i[f]'(z)= \displaystyle \frac{\partial Q_i[f]}{\partial x_0}(z)$ is a holomorphic map of complex variable $z=x_0+iy$ and $i\in \mathbb{S}.$
\end{definition}
\begin{remark}\label{eq:200}
Let $j\in \mathbb{S}$ be such that $j\bot i.$ Then there exist holomorphic functions  $f_1, f_2: \mathbb{B}_i \to \mathbb{C}(i)$ such that $Q_i[f]=f_1+f_2j $ and so $\displaystyle\frac{\partial f}{\partial x_0}(z)=f'_1(z)+f'_2(z)$ for some $z\in \mathbb{B}_i.$ Thus, for  $z\in \mathbb{B}_i,$  it follows that $$\left|f_l'(z)\right|^p\leq \left|\displaystyle\frac{\partial f}{\partial x_0}(z)\right|^p\leq 2^{{max}{(0,p-1)}}\left(|f'_1(z)|^p+|f'_2(z))|^p\right),~~l=1,2.$$ Thus, the function $f\in \mathfrak{B}_{p,i}$ if and only if $f_1,f_2 \in \mathfrak{B}_{p,\mathbb{C}}$ on $\mathbb{B}_i, $ (see \cite[Remark 4.3]{marco}). 
\end{remark}

\noindent The  proof of the  following proposition is analogus to  \cite[Proposition 2.6 ]{marco}.

\begin{proposition}\label{eq:20}
Let  $i\in  \mathbb{S},$  then $f\in  \mathfrak{B}_{p,i}, \; p>1$ if and only if $f\in  \mathfrak{B}_{p}.$ Moreover, the spaces $( \mathfrak{B}_{p,i},\|.\|_{ \mathfrak{B}_{p,i}})$ and $( \mathfrak{B}_p,\|.\|_{ \mathfrak{B}_p})$ have equivalent norms. More precisely, one has  $$\|f\|^p_{\mathfrak{B}_{p,i}}\leq \|f\|^p_{\mathfrak{B}_{p}}\leq 2^p\|f\|^p_{\mathfrak{B}_{p,i}}.$$
\end{proposition}
\noindent  For all $z,w\in \mathbb{D},$  Bergman metric on the unit disc  $\mathbb{D}$   in the complex plane $\mathbb{C}$ is given by  $$\beta(z,w)=\frac{1}{2}\log\frac{1+\rho(z,w)}{1-\rho(z,w)},$$ where $\rho(z,w)=|\frac{z-w}{1-\bar {z} w}|.$ 

\begin{definition}\cite{marco}.
  For $i\in \mathbb{S}$ and all $z,w\in \mathbb{B}_i,$ we define  $$\beta_i(z,w)=\displaystyle\frac{1}{2}\log\left(\frac{1+\frac{|z-w|}{|1-\bar z w|}}{1-\frac{|z-w|}{|1-\bar z w|}}\right).$$ 
\end{definition}

\begin{proposition}
For $1<  p,t< \infty,$ with $\frac{1}{p}+\frac{1}{t}=1.$ Let $f\in \mathfrak{B}_p$  and $i\in \mathbb{S}$ be fixed. Then for all  $q,w \in  \mathbb{B}_i,$ there exists  a constant $M_p>0$ such that $$\displaystyle|f(q)-f(w)|\leq 2M_p\|f\|_{\mathfrak{B}_{p}} \beta_i(q,w)^{\frac{1}{t}},$$  where $$\beta_i(q,w)=\displaystyle\frac{1}{2}\log\left(\frac{1+\frac{|q-w|}{|1-\bar q w|}}{1-\frac{|q-w|}{|1-\bar q w|}}\right). $$
\end{proposition}
\beginpf
 By Lemma  \ref{eq:1}, there exist two holomorphic functions $f_1, f_2 : \mathbb{B}_i \to \mathbb{C}(i)$ such that  $Q_i[f]=f_1+f_2j,$ where $j\bot i.$   Moreover,  the functions $f_l\in \mathfrak{B}_{p,\mathbb{C}}~;~l=1,2$ on $\mathbb{B}_i$. Furthermore, $\displaystyle\|f_{l}\|^{p}_{\mathfrak{B}_{p,\mathbb{C}}}\leq  \|f\|^{p}_{\mathfrak{B}_{p,i}}~;~ l=1,2~\mbox{and}~p>1.$  Therefore, from \cite[Theorem 9]{zhu91}, it follows that  for all  $q,w \in  \mathbb{B}_i$ in the  complex plane $\mathbb{C}(i),$ one has 
$$\begin{array}{ccl}
\displaystyle|f(q)-f(w)|^p&\leq&\displaystyle 2^{p-1}\left(|f_1(q)-f_1(w)|^p +|f_2(q)-f_2(w)|^p\right)\\
&\leq&\displaystyle 2^{p-1}M_p\left(\|f_1\|^{p}_{\mathfrak{B}_{p,\mathbb{C}}}\beta(q,w)^{\frac{p}{t}}+\|f_2\|^{p}_{\mathfrak{B}_{p,\mathbb{C}}}\beta(q,w)^{\frac{p}{t}}\right)\\
&\leq&\displaystyle  2^{p-1} 2M_p\|f\|^{p}_{\mathfrak{B}_{p,i}}\beta_i(q,w)^{\frac{p}{t}}\\
&\leq&\displaystyle  2^p M_p\|f\|^{p}_{\mathfrak{B}_{p}}\beta_i(q,w)^{\frac{p}{t}}.\\
\end{array}$$
\endpf
 \noindent  The following proposition on  Besov spaces over the unit disk was proved in   \cite[Theorem 8]{zhu91} and for its proof on Bloch spaces of slice holomorphic functions one can refer to   \cite[Theorem 2.20]{marco}.

\begin{proposition}
Let $f\in \mathfrak{B}_p, \; p >1$ and $\{a_n\}_{n\in \mathbb{N}}\subset \mathbb{H}$ be a sequence of quaternions such that $$f(q)=\sum_{n=0}^\infty q^ na_n  ~for~~q\in \mathbb{B}.$$ Then there exists a constant $K_p>0$ such that $$|a_n|^p\leq  \displaystyle2^p\frac{K_p}{n}\|f\|^p_{\mathfrak{B}_{p}}~~ for ~~n\in \mathbb{N}\cup \{0\}.$$
\end{proposition}
\beginpf
 Let $i,j\in \mathbb{S}$ be such that  $j\bot i.$  On applying Splitting Lemma 1.1, we can restrict $f$ on $\mathbb{B}_i$  such that $Q_i[f]=f_1+f_2j,$  for some  holomorphic functions $f_1, f_2 : \mathbb{B}_i \to \mathbb{C}(i)$ in the complex Besov space $\mathfrak{B}_{p, \mathbb{C}}$ on $\mathbb{B}_i.$ Furthermore, for any $z\in \mathbb{B}_i$ and $p>1,$ we have $$|f(z)|^p\leq 2^{p-1}(|f_1(z)|^p+|f_2(z)|^p).$$ Now for any $n\in \mathbb{N}\cup \{0\},$ let $a_{1,n}, \;a_{2,n}\in \mathbb{C}(i)$ such that $a_n=a_{1,n}+a_{2,n}j.$ Thus we have
$$|f(z)|^p=\displaystyle\left|\sum_{n=0}^\infty z^ na_n\right|^p\leq  2^{p-1}\left(\displaystyle\left|\sum_{n=0}^\infty z^ na_{1,n}\right|^p+\displaystyle\left|\sum_{n=0}^\infty z^ na_{2,n}\right|^p\right)=  2^{p-1}\left(|f_1(z)|^p+|f_2(z)|^p\right).$$ Therefore, from \cite[Theorem 8 (1)]{zhu91}, it follows that for any $n\in \mathbb{N},$ we have $$\displaystyle|a_{l,n}|\leq\frac{ K_p}{n^{\frac{1}{p}}}\|f_l\|_{\mathfrak{B}_{p,\mathbb{C}}}~~;~~l=1,2,$$ and so $\|f_l\|_{\mathfrak{B}_{p,\mathbb{C}}}\leq \|f\|_{\mathfrak{B}_{p,i}}~;~l=1,2.$ Then,    one has
$$\begin{array}{ccl}
|a_n|^p&=& \displaystyle2^{p-1}\left(|a_{1,n}|^p+|a_{2,n}|^p\right)\\
&\le& \displaystyle2^{p-1}\frac{K_p}{n}\left(\|f_1\|^p_{\mathfrak{B}_{p,\mathbb{C}}}+\|f_2\|^p_{\mathfrak{B}_{p,\mathbb{C}}}\right)\\
&\leq&  \displaystyle2^{p-1}2\frac{K_p}{n}\|f\|^p_{\mathfrak{B}_{p,i}}\\
&\leq&  \displaystyle2^p\frac{K_p}{n}\|f\|^p_{\mathfrak{B}_{p}}.\\
\end{array}$$
\endpf
\begin{remark}\label{eq:121}

Let  $L^p(\mathbb{B}_i,d\lambda_i, \mathbb{H }), \; 1\le p< \infty$   denote the space of quaternionic valued  equivalence classes of measurable  functions $g: \mathbb{B}_i\to \mathbb{H}$ such that  $$\int_{\mathbb{B}_i}|g(w)|^pd\lambda_i(w)< \infty.$$ Furthermore, for any $j\in \mathbb{S}$ with $j\bot i$ and $g=g_1+g_2j$  where $g_1,g_2 $ are holomorphic functions in complex plane $\mathbb{C}(i),$ then,  $g\in L^p(\mathbb{B}_i,d\lambda_i, \mathbb{H})$  if and only if $g_l\in L^p(\mathbb{B}_i,d\lambda_i, \mathbb{C}(i)), ~l=1,2,$  the usual $ L^p$-space of complex valued measurable functions on  $\mathbb{B}_i.$
\end{remark} 

  Now we define the bounded mean oscillation of the slice regular functions.
\begin{definition}
For any $z\in \mathbb{\mathbb{B}}_i,$ let $ \Delta_i(z,r)=\{w\in \mathbb{B}_i: \beta_i(z,w)<r\} \subset \mathbb{B}_i $ for some $r>0,$ be the Euclidean disk. Let $f^*_{r,i}(z)=\displaystyle\frac{1}{2\pi}\int_{ \Delta_i(z,r)}f(w)dA_i(w),$ for some arbitrary $i\in \mathbb{S}.$\\ 
A slice regular function $f$ is said to be in $BMO(\mathbb{B}_i)$ if $$\sup_{z \in \mathbb{B}_i}\frac{1}{2\pi}\displaystyle\int_{ \Delta_i(z,r)}|f(w)-f^*_{r,i}(z)|^pdA_i(w)<\infty,$$ with norm defined by $$\|f\|_{BMO(\mathbb{B}_i)}=\sup_{z \in \mathbb{B}_i}\left(\frac{1}{2\pi}\displaystyle\int_{\Delta_i(z,r)}|f(w)-f^*_{r,i}(z)|^pdA_i(w)\right)^{\frac{1}{p}}.$$ We say function $f\in BMO(\mathbb{B})$ if  $$\|f\|_{BMO(\mathbb{B})}:=\sup_{i\in  \mathbb{S}}\|f\|_{BMO(\mathbb{B}_i)}:=\sup_{i\in \mathbb{S}}\Lambda_{r,i}(f)<\infty,$$  where  $$\Lambda_{r,i}(f)(z)=\displaystyle\sup\{|f(z)-f(w)|: w\in  \Delta_i(z,r)\}~\mbox{for some}~ i\in  \mathbb{S}.$$
\end{definition}
\begin{proposition}
Let $p>1$  and $i,j \in  \mathbb{S}.$  Then $f\in BMO(\mathbb{B}_i)$ if and only if $f\in BMO(\mathbb{B}_j).$
\end{proposition}
\beginpf
Let $f\in SR(\mathbb{B})$ and choose $w=x+yj \in \mathbb{B}_j$ and $z=x+yi \in \mathbb{B}_i.$  Then by Representation formula, we have
 $$\left|f(w)\right| =\frac{1}{2}\left|(1-ji)f(z)+(1+ji)f(\bar z)\right|\leq \left|f(z)\right|+\left|f(\bar z)\right|.$$ Therefore
$$\begin{array}{ccl}
\displaystyle\frac{1}{2\pi}\int_{ \Delta_j(z,r)}|f(w)-f^*_{r,j}(z)|^pdA_j(w) &\leq &\displaystyle2^{max\{p-1,0\}}\frac{1}{2\pi}\displaystyle\int_{ \Delta_i(w,r)}|f(z)-f^*_{r,i}(w)|^pdA_i(z)\\
& +&\displaystyle 2^{max\{p-1,0\}}\frac{1}{2\pi}\displaystyle\int_{ \Delta_i(w,r)}|f(\bar z)-f^*_{r,i}( \bar w)|^pdA_i(\bar z).\\
\end{array}$$
On changing  $\bar z \to z$ and $\bar w \to w,$ we have
$$\displaystyle\frac{1}{2\pi}\int_{ \Delta_j(z,r)}\left|f(w)-f^*_{r,j}(z)\right|^pdA_j(w)\leq  2^{max\{p,1\}} \frac{1}{2\pi}\displaystyle\int_{ \Delta_i(w,r)}\left|f(z)-f^*_{r,i}(w)\right|^pdA_i(z).$$ Thus, we conclude that for any $f\in BMO(\mathbb{B}_i)$ implies $f\in BMO(\mathbb{B}_j).$ Finally, on interchanging the role of $i$ and $j$, we get the remaining one.
\endpf
\begin{proposition}
For $p>1$ and $\alpha>-1,$  let $f\in \mathfrak{B}_p .$ Then $f\in  BMO(\mathbb{B})$ if and only if  $ f\in BMO(\mathbb{B}_i),$ for some $i\in  \mathbb{S}.$
\end{proposition}
\beginpf 
Since the direct part is obivious, so we only remains to prove the converse part. suppose $f\in BMO(\mathbb{B}_i),$    for some arbitrary imaginary unit $i$ in $ \mathbb{S}.$ Therefore by  Representation formula, we have
$$\begin{array}{ccl}
\displaystyle\frac{1}{2\pi}\int_{ \Delta_j(z,r)}\left|f(w)-f^*_{r,j}(z)\right|^pdA_j(w)&\leq& 2^{p-1}\displaystyle\frac{1}{2\pi}\left(\int_{ \Delta_i(w,r)}\left|f(z)-f^*_{r,i}(w)\right|^pdA_i(z)\right)\\
&+&\displaystyle2^{p-1}\frac{1}{2\pi}\left(\int_{ \Delta_i(w,r)}\left|f(\bar z)-f^*_{r,i}(\bar w)\right|^pdA_i(\bar z)\right).
\end{array}$$
 On taking supremum over all $z\in \mathbb{B}_i,$ we have 
$$\begin{array}{ccl}
\displaystyle\|f\|_{BMO(\mathbb{B}_j)}&\leq&\displaystyle \sup_{z\in  \Delta_i(w,r)} 2^{p-1}\displaystyle\frac{1}{2\pi}\left(\int_{ \Delta_i(w,r)}\left|f(z)-f^*_{r,i}(w)\right|^pdA_i(z)\right)\\
&+&\displaystyle \displaystyle \sup_{z\in  \Delta_i(w,r)} 2^{p-1}\frac{1}{2\pi}\left(\int_{ \Delta_i(w,r)}\left|f(\bar z)-f^*_{r,i}(\bar w)\right|^pdA_i(\bar z)\right)\\
&\leq& 2^{p-1}2 \|f\|_{BMO(\mathbb{B}_i)}<\infty.
\end{array}$$

 Since  $j$ is  arbitrary, so we get   the desired result.\endpf

By previous proposition we conclude  the following inequality $$ \|f\|_{BMO(\mathbb{B}_i)}^p\leq  \|f\|_{BMO(\mathbb{B})}^p\leq 2^p \|f\|_{BMO(\mathbb{B}_i)}^p.$$
\begin{proposition}
For $p> 1,$  let $f$ be a slice regular function. Then $f\in \mathfrak{B}_p$ if and only if $\Lambda_r(f) \in L^p(\mathbb{B}_i,d\lambda_i, \mathbb{H }),$  for some $i\in \mathbb{S}.$
\end{proposition}
\beginpf
Suppose  $f\in \mathfrak{B}_p$ implies $f\in \mathfrak{B}_{p,i}.$ Let $j\in  \mathbb{S}$ be such that $j\bot i.$  By  Splitting Lemma (\ref{eq:1}), we can restrict $f$ on $\mathbb{B}_i$   with respect to  $j$, as $Q_i[f](z)=f_1(z)+f_2(z)j,$  for some holomorphic functions $f_1, f_2\in \mathbb{C}(i).$ If we decompose $\Lambda_r(f)$ on $\mathbb{B}_i$ as $\Lambda_r(f)=\Lambda_{r,1}(f_1)+\Lambda_{r,2}(f_2)j,$ for some complex  oscillation functions  $\Lambda_{r,1}(f_1)$  and $\Lambda_{r,2}(f_2)$. Then one can see directly from the complex result (see  \cite[Theorem 6]{zhu91} ) and Remark \ref{eq:121} that the functions $\Lambda_{r,l}(f_l) ;~ l=1,2$   lie in the usual $L^p$-space of complex valued measurable functions on $\mathbb{B}_i$ if and only if $\Lambda_r(f)\in L^p(\mathbb{B}_i,d\lambda_i, \mathbb{H }).$\\
 Conversely, assume $\Lambda_r(f)\in L^p(\mathbb{B}_i,d\lambda_i, \mathbb{H }).$  So we can write 
$$\begin{array}{ccl}
\Lambda_{r,1}(f_1)+\Lambda_{r,2}(f_2)j&=&\Lambda_r(f)\\
&=&\displaystyle\sup_{i\in  \mathbb{S}}\sup\{|f_1(z)-f_1(w)|: w\in  \Delta_i(z,r)\subset \mathbb{B}_i\}\\
&+&\displaystyle\sup_{i\in  \mathbb{S}}\sup\{|f_2(z)-f_2(w)|: w\in  \Delta_i(z,r)\subset \mathbb{B}_i\}.\\
\end{array}$$
This implies  $$\Lambda_{r,l}(f_l)=\displaystyle\sup_{i\in  \mathbb{S}}\sup\{|f_l(z)-f_l(w)|: w\in  \Delta_i(z,r)\subset \mathbb{B}_i\} \in L^p(\mathbb{B}_i,d\lambda_i, \mathbb{C }(i)),~~\mbox{for}~~l=1,2. $$ Again thanks to above classical result, we conclude that both $f_1$ and $f_2$ belong to complex Besov space $ \mathfrak{B}_{p, \mathbb{C }}$ on  $\mathbb{B}_i$ which is equivalent to $f\in \mathfrak{B}_{p,i}(\mathbb{B}_i)$ and so  $f\in \mathfrak{B}_p(\mathbb{B}).$
\endpf

\begin{proposition}\label{eq:501}
For $p>1,$ let $f\in SR(\mathbb{B}).$ Then $f\in \mathfrak{B}_p$ if and only if 
\beq\label{eq:207}
BMO(f) \in L^p(\mathbb{B}_i,d\lambda_i, \mathbb{H}),~ for~ some~ i\in \mathbb{S}.
\eeq
\end{proposition}
\beginpf
Let $f\in  \mathfrak{B}_p.$ Then $f\in \mathfrak{B}_{p,i}.$  Let $j\in \mathbb{S}$ with $j\bot i.$ According to Lemma (\ref{eq:1}) , any  $f\in SR(\mathbb{B})$ can  be restricted to  $\mathbb{B}_i$ decomposes  as $Q_i[f](z)=f_1(z)+f_2(z)j,$ for some $z\in \mathbb{B}_i$  and holomorphic functions $f_1,f_2 \in \mathbb{B}_i$. Thus, the  condition (\ref{eq:207}) holds if and only if $$BMO(f_l) \in L^p(\mathbb{B}_i,d\lambda_i, \mathbb{C}(i)),~ \mbox{for some}~ i\in \mathbb{S},~l=1,2.$$ Now, by  \cite[Theorem 7]{zhu91}, it follows that the above condition holds if and only if $f_1,f_2$ lie in the complex Besov space $\mathfrak{B}_{p, \mathbb{C}}$ on $\mathbb{B}_i$ which is same as $f\in \mathfrak{B}_{p,i}$ and so  $f\in \mathfrak{B}_{p}.$
\endpf

\section{Composition operators on Besov spaces}
\label{sec:3}
\subsection {Boundedness and Compactness}
\label{sec:3}
In this section, we characterize  boundedness and compactness of composition operators on Besov spaces of the slice holomorphic functions.
\begin{definition}
Let  $0<p<\infty $ and  let $\Phi: \mathbb{B}\rightarrow \mathbb{B}$ be a  slice hyperholomorphic  map such that $\Phi(\mathbb{B}_i)\subset \mathbb{B}_i$ for some $i \in \mathbb{S}.$ Then the composition operator
$C_\Phi$ on $\mathfrak{B}_p$ on the unit ball $\mathbb{B}$ induced by $\Phi$ is defined by $$C_\Phi f=f\circ_ i\Phi,~~for ~all~ f\in \mathfrak{B}_p.$$
\end{definition}
Composition operators are extensively studied on various holomorphic function spaces of different domains in $\mathbb{C}$ or $\mathbb{C}^n.$  For a study of composition operators on spaces of holomorphic functions, one can refer to  \cite{cowe} and \cite{shap1}.  For composiion operators on Besov spaces see,  \cite{arco02}.  A study of composition operators on Hardy spaces of slice holomorphic functions is initated  in \cite{grxw14}. Recently, Carleson measures for Hardy and Bergman spaces   in the quaternionic unit ball are characterized in \cite{saba16}. In \cite{arco14},  Hankel operators are studied on Hardy spaces via Carleson measures in a quaternionic variables. \\
\noindent The following theorem characterize  bounded composition operators  on the slice regular Besov spaces $\mathfrak{B}_{p}.$
\begin{theorem}
 Let $\Phi$ be  a slice holomorphic map on $\mathbb{B}$ such that $\Phi(\mathbb{B}_i)\subset \mathbb{B}_i$ for some $i\in \mathbb{S}.$  For all  $q\in \mathbb{B}$ and $a\in \mathbb{B}_i$,  let  $\sigma_a(q)=(1-qa)^*(a-q)$  be slice  regular M\"{o}bius transformation, Then the composition operator $C_{\Phi} $ is bounded on Besov space $ \mathfrak{B}_p, \; 1 <p<\infty$ if and only if  
\beq \label{eq:210}
\sup\|C_{\Phi}\sigma_a\|_{\mathfrak{B}_{p}}<\infty.
\eeq

\end{theorem}
\beginpf
Since  the  slice regular  M\"{o}bius transformation  on $\mathbb{B}_i$  coincides with the usual one dimensional complex  M\"{o}bius transformation, so assume    $\sigma_a\in \mathfrak{B}_{p,i}.$ Let $j\in \mathbb{S}$ with $j\bot i.$ So we can write  $\sigma_a=\sigma_{a,1}+\sigma_{a,2}j,$  for  each  one dimensional complex  M\"{o}bius transformation $\sigma_{a,l} \in  \mathfrak{B}_{p,  \mathbb{C}},~l=1,2.$

 Therefore, from \cite[ Theorem 13]{araz},  we have
\begin{eqnarray}\label{eq:213}
\displaystyle\sup_{i\in \mathbb{S}}\int_{\mathbb{B}_i}\left|(1-|z|^2)\frac{\partial C_\Phi \sigma_a}{\partial x_0}(z)\right|^pd\lambda_i(z)&\leq&2^{p-1}\displaystyle\sup_{i\in \mathbb{S}}\int_{\mathbb{B}_i}\left|(1-|z|^2)\frac{\partial C_\Phi \sigma_{a,1}}{\partial x_0}(z)\right|^pd\lambda_i(z) \nonumber\\ 
&+&2^{p-1}\displaystyle\sup_{i\in \mathbb{S}}\int_{\mathbb{B}_i}\left|(1-|z|^2)\frac{\partial C_\Phi \sigma_{a,2}}{\partial x_0}(z)\right|^pd\lambda_i(z) \nonumber  \nonumber\\
&=&\displaystyle2^{p-1}(\|C_{\Phi}\sigma_{a,1}\|_{\mathfrak{B}_{p,\mathbb{C}}}^p+\|C_{\Phi}\sigma_{a,2}\|_{\mathfrak{B}_{p,\mathbb{C}}}^p)  \nonumber\\
&\leq& 2^p\|C_{\Phi}\sigma_a\|_{\mathfrak{B}_{p,i}}^p. \\
 \nonumber\end{eqnarray}
Now, let $q=x_0+Iy \in \mathbb{B}$ for some $I\in \mathbb{S}.$ Then by  Theorem \ref{eq:115}, it follows that $$\left|\frac{\partial C_\Phi \sigma_a}{\partial x_0}(q)\right|=\left|\frac{1}{2}(1-Ii)\frac{\partial C_\Phi \sigma_a}{\partial x_0}(z)+\frac{1}{2}(1+Ii)\frac{\partial C_\Phi \sigma_a}{\partial x_0}(\bar z)\right|.$$ Since, as $|q|=|z|=|\bar z|,$ on applying triangle inequality, we have $$\left|(1-|q|^2)\frac{\partial C_\Phi \sigma_a}{\partial x_0}(q)\right|\leq\left|(1-|z|^2)\frac{\partial C_\Phi \sigma_a}{\partial x_0}(z)\right|+\left|(1-|\bar z|^2)\frac{\partial C_\Phi \sigma_a}{\partial x_0}(\bar z)\right|.$$ On taking integeral over $\mathbb{B}_i$ on both sides of  the above inequality  and for $p>1,$ we see
\begin{eqnarray}\label{eq:214}
\displaystyle\sup_{q\in \mathbb{B}}\int_{\mathbb{B}_i}\left|(1-|q|^2)\frac{\partial C_\Phi \sigma_a}{\partial x_0}(q)\right|^pd\lambda_i(q)&\leq&\displaystyle\sup_{i\in \mathbb{S}}\sup_{z\in \mathbb{B}_i}\int_{\mathbb{B}_i}\left|(1-|z|^2)\frac{\partial C_\Phi \sigma_a}{\partial x_0}(z)\right|^pd\lambda_i(z)\nonumber\\
&+&\displaystyle\sup_{i\in \mathbb{S}}\sup_{\bar z\in \mathbb{B}_i}\int_{\mathbb{B}_i}\left|(1-|\bar z|^2)\frac{\partial C_\Phi \sigma_a}{\partial x_0}(\bar z)\right|^pd\lambda_i(\bar z)\nonumber\\
&\leq& 2\displaystyle\sup_{i\in \mathbb{S}}\sup_{z\in \mathbb{B}_i}\int_{\mathbb{B}_i}\left|(1-|z|^2)\frac{\partial C_\Phi \sigma_a}{\partial x_0}(z)\right|^pd\lambda_i(z).
\end{eqnarray}
Thus, by using  (\ref{eq:213}) in  (\ref{eq:214}), we have 
$$\begin{array}{ccl}
\sup\|C_{\Phi}\sigma_a\|_{\mathfrak{B}_{p}}^p&=&\displaystyle\sup_{i\in \mathbb{S}}\sup_{q\in \mathbb{B}}\int_{\mathbb{B}_i}\left|(1-|q|^2)\frac{\partial C_\Phi \sigma_a}{\partial x_0}(q)\right|^pd\lambda_i(q)\nonumber\\
&\leq& 2\displaystyle\sup_{i\in \mathbb{S}}\sup_{z\in \mathbb{B}_i}\int_{\mathbb{B}_i}\left|(1-|z|^2)\frac{\partial C_\Phi  \sigma_a}{\partial x_0}(z)\right|^pd\lambda_i(z)\nonumber\\
&\leq&  2^{p+1} \|C_{\Phi} \sigma_a\|_{\mathfrak{B}_{p,i}}^p.\\ 
\end{array}$$
Since  $C_\Phi$ is bounded operator on the complex Besov space on $\mathbb{B}_i,$ therefore   $\|C_{\Phi}\sigma_a\|_{\mathfrak{B}_{p}}^p<\infty.$
Now suppose condition  (\ref{eq:210}) holds. Then by \cite[Theorem 13]{araz}, it holds if and only if $C_\Phi$ is bounded operator on the complex Besov space which is equivalent to the boundedness of $C_\Phi$ on $\mathfrak{B}_{p,i}$ and  so  $C_\Phi $ is bounded on   $\mathfrak{B}_p$.

 By using Splitting Lemma, Remark   \ref{eq:200} and  \cite[Lemma 3.6]{tjan}, the proof  of the following lemma follows easily.
\begin{lemma}\label{eq:208}
For $p\geq1,$ let $\mathfrak{B}_p$ be a slice regular Besov space on the unit ball $\mathbb{B}.$ Then the following condition holds:
\begin{itemize}
\item[(1)] every slice regular bounded sequence $\{f_m\}_{m\in \mathbb{N}}$ in $\mathfrak{B}_p$  on compact sets is uniformly bounded;
\item[(2)] for any  slice regular  sequence $\{f_m\}_{m\in \mathbb{N}}$ in $\mathfrak{B}_p$  such that $\|f_n\|_{\mathfrak{B}_p}\to 0,  \; f_n-f_n(0)\to 0$ uniformly on the compact sets.
\end{itemize}
\end{lemma}
 The next result is essential for the proof of   Proposition  \ref{eq:212}.
\begin{lemma}\label{eq:209}\cite[Lemma 3.7]{tjan} Let $X, Y$ be two Banach spaces of analytic functions on the unit disk $\mathbb{D}$. Suppose 
\begin{itemize}
 \item[(1)] the point evaluation functionals on X are continuous;
\item[(2)]  the closed unit ball in X is a compact subset of X in the topology of uniform convergence on  compact sets;
\item[(3)] $T: X\to Y$ is continuous,  where  X and Y  are equipped with  the topology of uniform convergence on compact sets.  
Then T is  a compact  operator if and only if given a bounded sequence $\{f_n\}$ in X such that $f_n \to 0$ uniformly on compact sets, then the sequence $\{Tf_n\}$ converges to zero in the norm of Y.
\end{itemize}
\end{lemma}
The   following proposition gives  the characterization  for compact composition operators.  
\begin{proposition}\label{eq:212}
For $p> 1,$ let $\mathfrak{B}_p$ be a slice regular Besov space on the unit ball $\mathbb{B}.$ Let $\Phi$ be  a slice holomorphic map on $\mathbb{B}$ such that $\Phi(\mathbb{B}_i)\subset \mathbb{B}_i$ for some $i\in \mathbb{S}.$ Then  $C_{\Phi}:\mathfrak{B}_p \to \mathfrak{B}_p$ is compact if and only if  for any  bounded sequence $\{f_m\}_{m\in \mathbb{N}}$ in $\mathfrak{B}_p $ with $f_m \to 0$ as $m\to \infty$ on compact sets, $\|C_{\Phi}f_m\| _{\mathfrak{B}_p }\to 0.$
\end{proposition}
\beginpf 
 The proof of the theorem is established if we prove the condtion of Lemma  \ref{eq:209}. As a consequence of Lemma \ref{eq:208}, we see that conditions (1) and (3) holds. Now, it  remains  to prove the  condition (2). For this, let $\{f_m\}$ be a slice regular bounded sequence in $\mathfrak{B}_p.$ Then by Lemma \ref{eq:208},  $\{f_m\}$ is uniformly bounded on the compact sets. Consider  $\{f_{m_k}\}$ a  subsequence of $\{f_m\}$ in $\mathfrak{B}_p$ such that  $f_{m_k}$ converges uniformly to $h$ on the compact sets, for some  $h\in SR(\mathbb{B}).$ Let $j \in \mathbb{S}$ with $j\bot i.$ Then by Lemma \ref{eq:1}, there exist holomorphic functions $f_{1,m_k},f_{2,m_k}: \mathbb{B}_i \to \mathbb{C}(i)$ such that $Q_i[f_{m_k}](z)=f_{1,m_k}(z)+f_{2,m_k}(z)j,$ for some $z\in \mathbb{B}_i.$  Furthermore, $f_{1,m_k}\to h_1$ and $f_{2,m_k}\to h_2$  uniformly on the compact sets, where $h_l\in \mathbb{C}(i), l=1,2$ with $Q_i[h]=h_1+h_2j.$ From Remark \ref{eq:200}, we  conclude that  $f_{1,m_k}$ and $f_{1,m_k}$  belong to the complex Besov space $ \mathfrak{B}_{p,  \mathbb{C}}( \mathbb{B}_i).$ Thus, from  \cite[Lemma 3.8]{tjan} and by applying Minkowski's inequality and Fatou's Theorem, for $p> 1$, we have 
$$\begin{array}{ccl}
\displaystyle\left(\int_{\mathbb{B}_i}\left|\Bigg(\frac{\partial h}{\partial x_0}(z)\Bigg)(1-|z|^2)\right|^p d\lambda_i(z)\right)^\frac{1}{p}&\leq&\displaystyle\left(\int_{\mathbb{B}_i}2^{p-1}\left|(h'_1(z)+h'_2(z)j)(1-|z|^2)\right|^p d\lambda_i(z)\right)^\frac{1}{p}\\
&\leq&\displaystyle\left(2^{p-1}\int_{\mathbb{B}_i}\left|(h'_1(z))(1-|z|^2)\right|^p d\lambda_i(z)\right)^\frac{1}{p}\\
&+&\displaystyle\left(2^{p-1}\int_{\mathbb{B}_i}\left|(h'_2(z))(1-|z|^2)\right|^p d\lambda_i(z)\right)^\frac{1}{p}\\
&=&\displaystyle2^\frac{p-1}{p}\left(\int_{\mathbb{B}_i}\lim_{k\to \infty}\left|f_{1,m_k}'(z)(1-|z|^2)\right|^p d\lambda_i(z)\right)^\frac{1}{p}\\
&+&\displaystyle2^\frac{p-1}{p}\left(\int_{\mathbb{B}_i}\lim_{k\to \infty}\left|f_{2,m_k}'(z)(1-|z|^2)\right|^p d\lambda_i(z)\right)^\frac{1}{p}\\
&\leq&\displaystyle2^\frac{p-1}{p}\lim_{k\to \infty}\left(\int_{\mathbb{B}_i}\left|f_{1,m_k}'(z)(1-|z|^2)\right|^p d\lambda_i(z)\right)^\frac{1}{p}\\
&+&\displaystyle2^\frac{p-1}{p}\lim_{k\to \infty}\left(\int_{\mathbb{B}_i}\left|f_{2,m_k}'(z)(1-|z|^2)\right|^p d\lambda_i(z)\right)^\frac{1}{p}\\
&=&\displaystyle 2^\frac{p-1}{p}\left(\lim_{k\to \infty}\inf \|f_{1,m_k}\|_{\mathfrak{B}_{p,  \mathbb{C}}}+\lim_{k\to \infty}\inf\|f_{2,m_k}\|_{\mathfrak{B}_{p,  \mathbb{C}}}\right)\\
&\leq& 2^\frac{2p-1}{p}\displaystyle\lim_{k\to \infty}\inf\left(\|f_{m_k}\|_{\mathfrak{B}_{p,i}}\right)\\
&<& \infty.
\end{array}$$
\endpf
\noindent The next result is the immediate consequence of Proposition \ref{eq:212}.
\begin{corollary}
For $1< p< \infty,$  let $\Phi$ be a slice holomorphic map  such that $\Phi(\mathbb{B}_i)\subset \mathbb{B}_i$ for some $i\in \mathbb{S}.$  If $\|\Phi\|_{\infty}<1,$ then  $C_{\Phi}:\mathfrak{B}_p \to \mathfrak{B}_p$ is  comapct. 
\end{corollary}
\beginpf
Let $\{f_n\}$ be a bounded sequence in $\mathfrak{B}_p .$ Then $f_n \in \mathfrak{B}_{p,i} $ such that $f_n \to 0$ uniformly on the compact  subsets of $ \mathbb{B}_i$ for some $i \in \mathbb{S}.$  Let $j\in \mathbb{S}$ be such that $j\bot i.$  Let $f_{1,n},f_{2,n}: \mathbb{B}_i \to \mathbb{C}(i)$   be holomorphic functions  
such that  $Q_i[f](z)=f_{1,n}(z)+f_{2,n}(z)j,$ for some  $z=x_0+iy\in \mathbb{B}_i.$ By   Remark \ref{eq:200}, we have   $f_{1,n},  \;f_{2,n}$ lie in the complex  Besov space $\mathfrak{B}_{p,\mathbb{C}}$ on $\mathbb{B}_i,$ where $\mathbb{B}_i$ is identified with $\mathbb{D}\subset \mathbb{C}(i).$ Therefore,
\begin{eqnarray}\label{eq:230}
\displaystyle\sup_{i\in \mathbb{S}}\int_{\mathbb{B}_i}\left|(1-|z|^2)\frac{\partial C_\Phi f_n}{\partial x_0}(z)\right|^pd\lambda_i(z)&\leq&2^{p-1}\displaystyle\sup_{i\in \mathbb{S}}\int_{\mathbb{B}_i}\left|(1-|z|^2)\frac{\partial C_\Phi f_{1,n}}{\partial x_0}(z)\right|^pd\lambda_i(z) \nonumber\\ 
&+&2^{p-1}\displaystyle\sup_{i\in \mathbb{S}}\int_{\mathbb{B}_i}\left|(1-|z|^2)\frac{\partial C_\Phi f_{2,n}}{\partial x_0}(z)\right|^pd\lambda_i(z) \nonumber  \nonumber\\
&=&\displaystyle2^{p-1}(\|C_{\Phi}f_{1,n}\|_{\mathfrak{B}_{p,\mathbb{C}}}^p+\|C_{\Phi}f_{2,n}\|_{\mathfrak{B}_{p,\mathbb{C}}}^p)  \nonumber\\
&\leq& 2^p\|C_{\Phi}f_{n}\|_{\mathfrak{B}_{p,i}}^p. \\
 \nonumber\end{eqnarray}
Now, appealing to   Theorem \ref{eq:115} and the fact that $|q|=|\bar z|= |z|,$   equation (\ref{eq:230}) and  \cite[Corollary 2.12]{tjan1996}, it follows that \\
$$\begin{array}{ccl}
\displaystyle\sup_{q\in \mathbb{B}}\int_{\mathbb{B}_i}\left|(1-|q|^2)\frac{\partial C_\Phi f_n}{\partial x_0}(q)\right|^pd\lambda_i(q)&\leq&\displaystyle\sup_{i\in \mathbb{S}}\sup_{z\in \mathbb{B}_i}\int_{\mathbb{B}_i}\left|(1-|z|^2)\frac{\partial C_\Phi f_n}{\partial x_0}(z)\right|^pd\lambda_i(z)\\
&+& \displaystyle\sup_{i\in \mathbb{S}}\sup_{\bar z\in \mathbb{B}_i}\int_{\mathbb{B}_i}\left|(1-|\bar z|^2)\frac{\partial C_\Phi f_n}{\partial x_0}(\bar z)\right|^pd\lambda_i(\bar z)\\
&\leq &\displaystyle 2^p\displaystyle\sup_{i\in \mathbb{S}}  \sup_{z\in \mathbb{B}_i}\int_{\mathbb{B}_i} \left|(1-|z|^2)\frac{\partial C_\Phi f_n}{\partial x_0}(z)\right|^pd\lambda_i(z)\\
&\leq & 2^{p+1}\displaystyle\sup_{i\in \mathbb{S}} \sup_{z\in \mathbb{B}_i}\int_{\mathbb{B}_i}\left|\frac{\partial  f_n}{\partial x_0}(\Phi(z))\right|^p\left|(1-|z|^2)\right|^p\\
&.&\displaystyle\left|\frac{\partial \Phi }{\partial x_0}(z)\right|^pd\lambda_i(z)\\
&\leq& 2^{p+1}\epsilon.
\end{array}$$
Therefore, $\|C_{\Phi}f_{m}\|_{\mathfrak{B}_{p}}^p \to 0$ as $n \to \infty.$  Hence the result.
\endpf
\noindent  The following proposition gives the compactness between  Besov and Bloch spaces of slice regular functions.
\begin{proposition}
 For $p>1,$ let $\Phi$ be  a slice holomorphic map on $\mathbb{B}$ such that $\Phi(\mathbb{B}_i)\subset \mathbb{B}_i,$ for some $i\in \mathbb{S}.$ Then $C_{\Phi}:\mathfrak{B}_p \to \mathcal{B}$ is  compact  if and only if  
\beq \label{eq:211}
\|C_{\Phi}\sigma_a\|_{\mathcal{B}}\to 0, ~~as ~~|a|\to 1,
\eeq
where    $\sigma_a(q)=(1-q\bar a)^**(a-q),~~q\in  \mathbb{B}$ and $\mathcal{B}$ is a slice regular Bloch space on the unit ball $\mathbb{B}.$ Here $\star$ denotes the slice regular product. 
\end{proposition}
\beginpf
 Let $\{\sigma_a:a\in \mathbb{B}\}$ be a set  in $\mathfrak{B}_{p}$ such that $\sigma_a-a\to 0$ as $|a| \to 1$. Suppose $ C_{\Phi}$ is compact operator. Then by Lemma \ref{eq:212}, $\{\sigma_a\}$ is a bounded set in $\mathfrak{B}_{p}.$  Therefore, 
$\|C_{\Phi}\sigma_a\|_{\mathcal{B}}= 0$. Suppose condition (\ref{eq:211}) holds. Let ${f_m}$ be   a bounded sequence in   $\mathfrak{B}_{p,i}$ such that $f_m\to 0$ uniformly on the compact sets  as $ m\to \infty.$  We  claim  $C_{\Phi}:\mathfrak{B}_p \to \mathcal{B}$ is  compact. For this, take $j\in \mathbb{S}$ with $j\bot i.$ Let $f_{1,m},f_{2,m}$ be holomorphic functions such that $Q_i[f_m]=f_{1,m}(z)+f_{2,m}(z)j,$ for some  $z=x_0+iy\in \mathbb{B}_i.$ By   Remark \ref{eq:200}, we have  $f_{1,m}, f_{2,m} $ lie in the complex Besov space $ \mathfrak{B}_{p,  \mathbb{C}}(\mathbb{B}_i).$ Therefore, from \cite[Theorem 4.1]{tjan} and as $\|f_{l}\|_{\mathcal{B}_{p,\mathbb{C}}}\leq \|f\|_ {\mathcal{B}_{p,i}},$ we have
$$\begin{array}{ccl}
\|C_{\Phi}f_m\|_{\mathcal{B}}&=&\displaystyle \sup_{i\in \mathbb{S}}\sup_{z\in \mathbb{B}_i}\left\{(1-|z|^2)\left|\frac{ \partial C_{\Phi}(f_{1,m}+f_{2,m}j)}{\partial x_0}(z)\right|\right\}\\
&=&\displaystyle \sup_{i\in \mathbb{S}} \sup_{z\in \mathbb{B}_i}\left\{(1-|z|^2)\left|\frac{ \partial C_{\Phi}f_{1,m}}{\partial x_0}(z)+\frac{ \partial C_{\Phi}f_{2,m}}{\partial x_0}(z)j\right| \right\}\\
&\leq&\displaystyle \sup_{i\in \mathbb{S}} \sup_{z\in \mathbb{B}_i}\left\{(1-|z|^2)\left|\frac{ \partial C_{\Phi}f_{1,m}}{\partial x_0}(z)\right|\right\}+\displaystyle  \sup_{i\in \mathbb{S}}\sup_{z\in \mathbb{B}_i}\left\{(1-|z|^2)\left|\frac{ \partial C_{\Phi}f_{2,m}}{\partial x_0}(z)\right|\right\}\\
&\leq&2\displaystyle \sup_{i\in \mathbb{S}}\sup_{z\in \mathbb{B}_i}\left\{(1-|z|^2)\left|\frac{ \partial C_{\Phi}f_{m}}{\partial x_0}(z)\right|\right\}\\
&=&2\displaystyle \sup_{i\in \mathbb{S}}\left\{\frac{(1-|z|^2)}{(1-|\Phi(z)|^2)}\left|\frac{ \partial \Phi}{\partial x_0}(z)\right| \sup_{z\in \mathbb{B}_i}(1-|\Phi(z)|^2)\left|\frac{ \partial f_m}{\partial x_0}(\Phi(z))\right|\right\}\\

&\leq&2\displaystyle \sup_{i\in \mathbb{S}}\left\{\frac{(1-|z|^2)}{(1-|\Phi(z)|^2)}\left|\frac{ \partial \Phi}{\partial x_0}(z)\right|\right\}\|f_m\|_{\mathcal{B}_i}\\
&\leq&2\displaystyle \sup_{i\in \mathbb{S}}\left\{ \frac{(1-|z|^2)}{(1-|\Phi(z)|^2)}\left|\frac{ \partial \Phi}{\partial x_0}(z)\right|\right\}\|f_m\|_{\mathfrak{B}_{p,i}}.\\
\end{array}$$
Since $\{f_m\}$ is bounded  in $\mathfrak{B}_{p,i},$ so $\|C_{\Phi}f_m\|_{\mathfrak{B}_{p,i}}\to 0$ as $m\to \infty.$ Thus, $\|C_{\Phi}f_m\|_{\mathcal{B}}\to 0$ as $m\to \infty.$ Hence by Lemma \ref{eq:212}, $C_{\Phi}:\mathfrak{B}_p \to \mathcal{B}$ is  compact.
\endpf
\section{Essential norm}
In this section, we find some estimates for the essential norm of composition operators on the slice regular Besov space. Firstly, we define Carelson measure. 
\begin{definition}
For $1<p<\infty,$ let $\mathfrak{B}_p $  be a slice regular Besov space. Let  $\mu$ be a  $\mathbb{H}$-valued  positive measure on $\mathbb{B}_i.$ Then  $\mu$ is said to be $\mathbb{H}$-valued  p-Carleson measure on $\mathbb{B}$ if there is a constant $M>0$ such that $$\int_{\mathbb{B}_i}\left|\frac{\partial f}{\partial x_0}(q)\right|^pd\mu(q)\leq M\|f\|_{\mathfrak{B}_p}^p,$$ for all $f\in \mathfrak{B}_p(\mathbb{B}).$
\end{definition}
\begin{theorem}
 Let $f\in SR(\mathbb{B}). $ If $\mu=\mu_1+\mu_2j$ for some $i\in \mathbb{S}.$ Then $\mu$ is $\mathbb{H}$-valued  p-carleson measure  on the slice regular Besov space if and only if $\mu_1,\mu_2$ are  p-Carleson measure  on the complex Besov space $\mathfrak{B}_{p.\mathbb{C}}, \;1<p<\infty $ in $\mathbb{B}_i.$
\end{theorem}
\beginpf
Let $j\in \mathbb{S}$ be such that $i\bot j.$ Then for any $f\in \mathfrak{B}_{p,i}$ there exist holomorphic functions $f_1,f_2: \mathbb{B}_i \to \mathbb{C}(i)$ 
such that  $Q_i[f]=f_{1}(z)+f_{2}(z)j,$ for some  $z=x_0+iy\in \mathbb{B}_i.$ Now,  $\mu$ is $\mathbb{H}$-valued  p-Carleson measure on $\mathfrak{B}_{p}$ if and only if  $\mu$ is $\mathbb{H}$-valued  p-Carleson measure on $\mathfrak{B}_{p,i}$  if and only if 
$$\begin{array}{ccl}
 \displaystyle\int_{\mathbb{B}_i}\left|\frac{\partial f}{\partial x_0}(q)\right|^pd\mu(q)&\leq& M\|f\|_{\mathfrak{B}_{p,i}}^p\\
&\Leftrightarrow&\displaystyle\int_{\mathbb{B}_i}\left|\frac{\partial( f_1+f_2j)}{\partial x_0}(q)\right|^pd(\mu_1+\mu_2j)(q)\\
&\leq& M\|f_1+f_2j\|_{\mathfrak{B}_{p,i}}^p\\
&\Leftrightarrow&\displaystyle\int_{\mathbb{B}_i}\left|\frac{\partial f_l}{\partial x_0}(q)\right|^pd\mu_l(q)\\
&\leq &2^p M\|f_l\|_{\mathfrak{B}_{p,\mathbb{C}}}^p
\end{array}$$
if and only if $\mu_l,$ for $l=1,2$ is p-Carleson measure on $\mathfrak{B}_{p,\mathbb{C}}(\mathbb{B}_i).$
\endpf
\begin{definition} \cite{fcs13, marco}
The slice regular M\"{o}bius transformation $\sigma_a$  for every $a\in \mathbb{B} $ is define  as  $$\sigma_a(q)=(1-qa)^{-*}*(a-q),~~for~q\in  \mathbb{B},$$ where $*$  is slice regular product.
\end{definition}
\noindent The  slice regular M\"{o}bius transformation $\sigma_a$  satisfies the following conditions:
\begin{itemize}
\item[(i)]  $\sigma_a:  \mathbb{B}\to\mathbb{B}$ is a bijective mapping,
\item[(ii)]  For all $z\in \mathbb{B}_i,~\sigma_a(z)=\displaystyle\frac{a-z}{1-\bar a z},$ 
\item[(iii)]  For all $q\in \mathbb{B},~ \sigma_a(a)=0, \sigma_a(0)=a$ and   $\sigma_a \circ \sigma_a(q)=q.$
\end{itemize}
Now we give the definition of essential norm.\
\begin{definition}
The essential norm of a continuous linear operator T between the normed linear spaces X and Y is its distance from the compact operator K, that is $$\|T\|_e^{X\to Y}=\inf \{\|T-K\|^{X\to Y}: K~\mbox{ is compact operator}\},$$ where $\|.\|^{X\to Y}$ denotes the operator norm and $\|.\|_e^{X\to Y}$ is the essential norm.
\end{definition}
\noindent  Here, we give an  essential norm estimate for composition operators  on the slice regular Besov space $\mathfrak{B}_{p}.$
\begin{theorem}
For $1<p<\infty$ and $\alpha>-1,$ let $\Phi$ be a slice holomorphic map  such that $\Phi(\mathbb{B}_i)\subset \mathbb{B}_i,$ for some $i\in \mathbb{S}.$ Suppose  the composition operator $C_{\Phi}:\mathfrak{B}_p \to \mathfrak{B}_p$ is bounded.
 Then there is an absolute  constant $M\geq 1$ such that  $$\lim_{|a|\to 1}\sup \sup_{i\in \mathbb{S}}\int_{\mathbb{B}_i}\frac{(1-|a|^2)^{2+\alpha}}{|1-\bar a q|^{2(2+\alpha)}}d\mu_{p}(q)\leq \| C_{\Phi}\|_{e}\leq M2^{p}\lim_{|a|\to 1}\sup \sup_{i\in \mathbb{S}}\int_{\mathbb{B}_i}\frac{(1-|a|^2)^{2+\alpha}}{|1-\bar a q|^{2(2+\alpha)}}d\mu_{p}(q).$$
\end{theorem}
\beginpf
Let  $f=\displaystyle \sum_{k=0}^\infty q^ka_k \in  \mathfrak{B}_{p,i},$ for some $i\in \mathbb{S}$. For $0<r<1,$  denote $\mathbb{B}_r=\{|z|<r\}$ in the complex plane $\mathbb{C}(i).$ Consider an operator $R_nf(q)=\displaystyle \sum_{k=n+1}^\infty q^ka_k,$ for some integer $n.$ Suppose   $j\in \mathbb{S}$ with  $j\bot i.$ Then there exists holomorphic functions $f_1,f_2: \mathbb{B}_i \to \mathbb{C}(i)$ 
such that  $Q_i[f]=f_{1}(z)+f_{2}(z)j,$ for some  $z=x_0+iy\in \mathbb{B}_i.$ By   Remark \ref{eq:200}, we have  $f_l =\displaystyle \sum_{k=0}^\infty q^ka_{l,k}\in \mathfrak{B}_{p,  \mathbb{C}}(\mathbb{B}_i),$ and $R_{l,n}f_l(q)=\displaystyle \sum_{k=n+1}^\infty q^ka_{l,k},$ for some integer $n$ and $ l=1,2.$ Therefore, we have  
$$\begin{array}{ccl}
\| C_{\Phi}\|_{e}&\leq& \displaystyle\lim_{n\to \infty }\inf\| C_{\Phi}R_n\|_{\mathfrak{B}_{p}}^p \leq \lim_{n\to \infty }\inf \sup_{\|f\|_{\mathfrak{B}_{p}}\leq 1}\| (C_{\Phi}R_n)f\|_{\mathfrak{B}_{p}}^p.
\end{array}$$
Now, 
$$\begin{array}{ccl}
\| (C_{\Phi}R_n)f\|_{\mathfrak{B}_{p}}^p&=&\left(\displaystyle\left|R_{1,n}f_1(\Phi(0))\right|+\sup_{i\in \mathbb{S}}\int_{ \mathbb{B}_i}\left|(1-|z|^2)\frac{\partial (C_{\Phi}R_{1,n})f_1}{\partial x_0}(\Phi(z))\right|^pd\lambda_i(z)\right)\\
&+&\left( \displaystyle|R_{2,n}f_2(\Phi(0))| +\sup_{i\in \mathbb{S}}\int_{ \mathbb{B}_i}\left|(1-|z|^2)\frac{\partial (C_{\Phi}R_{2,n})f_2}{\partial x_0}(\Phi(z))\right|^pd\lambda_i(z)j\right).
\end{array}$$ 
 Let $\mu=\mu_{1,p}+\mu_{2,p} j ,$ where $\mu_{1,p}$ and  $\mu_{2,p}$ are  two   p-Carleson measure on $ \mathbb{B}_i$ with the values in $ \mathbb{C}(i).$ Then again thanks to Remark \ref{eq:200} and \cite[Theorem 3.4]{sk09}, we have\\
$$\begin{array}{ccl}
\displaystyle\sup_{i\in  \mathbb{S}}\int_{ \mathbb{B}_i}\left|(1-|z|^2)\frac{\partial (C_{\Phi}R_{n})f}{\partial x_0}(\Phi(z))\right|^pd\lambda_i(z)&\leq&\displaystyle 2^{p-1}\sup_{i\in  \mathbb{S}}\int_{ \mathbb{B}_i}\left|(1-|z|^2)\frac{\partial (C_{\Phi}R_{1,n})f_1}{\partial x_0}(\Phi(z))\right|^pd\lambda_i(z)\\ 
\end{array}$$ 
$$\begin{array}{ccl}
&+&\displaystyle 2^{p-1}\sup_{i\in  \mathbb{S}}\int_{ \mathbb{B}_i}\left|(1-|z|^2)\frac{\partial (C_{\Phi}R_{2,n})f_2}{\partial x_0}(\Phi(z))\right|^pd\lambda_i(z).\\ 
&=&\displaystyle 2^{p-1}\sup_{i\in  \mathbb{S}}\int_{ \mathbb{B}_i}\left|\frac{\partial (C_{\Phi}R_{1,n})f_1}{\partial x_0}(q)\right|^pd\mu_{1,p}(z)\\ 
&+&\displaystyle 2^{p-1}\sup_{i\in  \mathbb{S}}\int_{ \mathbb{B}_i}\left|\frac{\partial (C_{\Phi}R_{2,n})f_2}{\partial x_0}(q)\right|^pd\mu_{2,p}(z).\\ 
&=& \displaystyle2^{p-1}\left(\| (C_{\Phi}R_{1,n})f_1\|_{\mathfrak{B}_{p, \mathbb{C}}}^p+\| (C_{\Phi}R_{2,n})f_2\|_{\mathfrak{B}_{p, \mathbb{C}}}^p \right)\\
&\leq&\displaystyle 2^{p} \| (C_{\Phi}R_n)f\|_{\mathfrak{B}_{p,i}}^p\\
&=&\displaystyle2^{p}\int_{ \mathbb{B}_i \setminus  \mathbb{B}_r}\left|\frac{\partial (C_{\Phi}R_{n})f}{\partial x_0}(q)\right|^pd\mu_{p}(z)\\
&+&\displaystyle2^{p}\int_{  \mathbb{B}_r}\left|\frac{\partial (C_{\Phi}R_{n})f}{\partial x_0}(q)\right|^pd\mu_{p}(z).
\end{array}$$
 Again by    \cite[Theorem 3.4]{sk09}, for some fixed $r,$ we have$$\displaystyle2^{p}\sup\int_{  \mathbb{B}_r}\left|\frac{\partial (C_{\Phi}R_{n})f}{\partial x_0}(q)\right|^pd\mu_{p}(z)\to 0~~\mbox{as}~~n\to \infty$$ and $$\displaystyle\int_{ \mathbb{B}_i \setminus  \mathbb{B}_r}\left|\frac{\partial (C_{\Phi}R_{n})f}{\partial x_0}(q)\right|^pd\mu_{p,r}(z)\leq C_1C_2\|\mu_{p}\|^*_{r},$$ for some absolute constants $C_1,C_2 $ and $\|\mu_{p}\|^*_{r}=\displaystyle\lim_{a\to 1}\sup\int_{ \mathbb{B}_i}\left|\frac{\partial \sigma_a(q) }{\partial x_0}\right|^pd\mu_p(q), $ for $0<r<1$ and $p>1.$  Let $\mu_{p,r}$ be the  restriction  of measure  $\mu_{p}$ to the set $ \mathbb{B}_i \setminus  \mathbb{B}_r.$ Thus, 
$$\begin{array}{ccl}
\|C_\Phi\|_e^p &\leq&  \displaystyle\lim_{n\to \infty }\inf \sup_{\|f\|_{\mathfrak{B}_{p}}\leq 1}\| (C_{\Phi}R_n)f\|_{\mathfrak{B}_{p}}^p\\
&\leq& 2^pC_1C_2\displaystyle\lim_{r\to 1}\|\mu_{p}\|^*_{r}\\
&=&2^pC_1C_2\displaystyle\lim_{a\to 1}\sup_{i\in  \mathbb{S}}\int_{ \mathbb{B}_i}\left|\frac{\partial \sigma_a(q) }{\partial x_0}\right|^pd\mu_p(q)\\
&=&\displaystyle 2^p\lim_{|a|\to 1}\sup_{i\in  \mathbb{S}}\int_{ \mathbb{B}_i}\frac{(1-|a|^2)^{2+\alpha}}{|1-\bar a q|^{2(2+\alpha)}}d\mu_{p}(q).
\end{array}$$
Thus, we have the upper bound.\\
Now, let $\sigma_a(z)=\displaystyle\frac{a-z}{1-\bar a z}$ be the complex M\"{o}bius transformation on $ \mathbb{B}_i,$ associated with $a$. Clearly $\sigma_a$ is bounded in $\mathfrak{B}_{p,i}.$ Also $\sigma_a-a \to 0$ as $|a|\to 1$ uniformly on the compact subsets of $ \mathbb{B}_i$ and $ |\sigma_a(z)-a|=\displaystyle |z|\frac{1-|a|^2}{|1-\bar a z|}.$ Furthermore, $\|K(\sigma_a-a)\|_{\mathfrak{B}_{p,i}}\to 0$ as $|a|\to 1$ for some compact operator K on $\mathfrak{B}_{p,i}.$ Therefore,  
$$\begin{array}{ccl}
\|C_\Phi\|_e^p&\geq& \|C_\Phi -K\|_{\mathfrak{B}_{p}}^p\geq \|C_\Phi -K\|_{\mathfrak{B}_{p,i}}^p\\
&\geq&  \displaystyle\lim_{|a|\to 1}\|(C_\Phi -K)\sigma_a\|_{\mathfrak{B}_{p,i}}^p\\
&\geq & \displaystyle\lim_{|a|\to 1}\sup\|C_\Phi \sigma_a\|_{\mathfrak{B}_{p,i}}^p-\lim_{|a|\to 1}\sup\|K\sigma_a\|_{\mathfrak{B}_{p,i}}^p\\
&=&\displaystyle\lim_{|a|\to 1}\sup \sup_{i\in \mathbb{S}}\int_{\mathbb{B}_i}\frac{(1-|a|^2)^{2+\alpha}}{|1-\bar a q|^{2(2+\alpha)}}d\mu_{p}(q).
\end{array}$$
Hence the desired result.
\endpf
\begin{remark}
{\normalfont {By using Splitting Lemma $1.1$ and  Representation Theorem $1.2$  it can be proved that the composition operators on spaces of  slice holomorphic will be bounded if and only if the corresponding composition operators are bounded on classical holomorphic function spaces. }}
\end{remark}

\end{document}